\documentclass{amsproc}
\usepackage{amsfonts}
\usepackage{graphicx}
\usepackage{amscd}

\theoremstyle{plain}

\newtheorem{remark}{Remark}

\newtheorem{theorem}{Theorem}
\numberwithin{equation}{section}

\input{tcilatex}

\begin{document}
\title[Reproducing kernels ]{The reproducing kernel structure arising from a
combination of continuous and discrete orthogonal polynomials into Fourier
systems}
\author{Lu\'{\i}s Daniel Abreu}
\address{Department of Mathematics, Universidade de Coimbra\\
Portugal}
\email{daniel@mat.uc.pt.}
\thanks{
Partial financial assistance by FCG, FCT post-doctoral grant
SFRH/BPD/26078/2005 and CMUC.}
\dedicatory{In memory of Joaquin Bustoz}

\begin{abstract}
We study mapping properties of operators with kernels defined via a
combination of continuous and discrete orthogonal polynomials, which provide
an abstract formulation of quantum ($q$-) Fourier type systems. We prove
Ismail%
\'{}%
s conjecture regarding the existence of a reproducing kernel structure
behind these kernels, by establishing a link with Saitoh%
\'{}%
s theory of linear transformations in Hilbert space. The results are
illustrated with Fourier kernels with ultraspherical weights, their
continuous $q$-extensions and generalizations. As a byproduct of this
approach, a new class of sampling theorems is obtained, as well as Neumann
type expansions in Bessel and $q$-Bessel functions.
\end{abstract}

\maketitle

\section{Introduction}

The Gegenbauer expansion of the two variable complex exponential in terms of
the ultraspherical polynomials 
\begin{equation}
e^{ixt}=\Gamma (\nu )\left( \frac{t}{2}\right) ^{-\nu }\sum_{k=0}^{\infty
}i^{k}(\nu +k)J_{\nu +k}(t)C_{k}^{\nu }(x)  \label{Geg}
\end{equation}%
has the remarkable feature of being at the same time an expansion in a
Neumann series of Bessel functions. The usefulness of this expansion was
made very clear in a paper authored by Ismail and Zhang, where it was used
to solve the eigenvalue problem for the left inverse of the differential
operator, on $L^{2}$ spaces with ultraspherical weights \cite{IZ}. The
consideration of the $q$-analogue of this diagonalization problem led the
authors to extend Gegenbauer%
\'{}%
s formula to the $q$-case. This task required the introduction of a new $q$%
-analogue of the exponential, a two variable function denoted by $\mathcal{E}%
_{q}(x;t)$ which became known in the literature as the curly $q$-exponential
function, bearing the name from its notational convention. Ismail and Zhang%
\'{}%
s formula is 
\begin{equation}
\mathcal{E}_{q}(x;it)=\frac{t^{-\nu }(q;q)_{\infty }}{(-qt^{2};q^{2})_{%
\infty }(q^{\nu +1};q)_{\infty }}\sum_{k=0}^{\infty }i^{k}q^{k^{2}/4}\frac{%
(1-q^{k+\nu })}{(1-q^{\nu })}J_{\nu +k}^{(2)}(2t;q)C_{k}(x;q^{\nu }|q).
\label{qGeg}
\end{equation}%
The functions involved in this formula will be defined in section 4. Since
its introduction, the function $\mathcal{E}_{q}$ was welcomed as a proper $q$%
-analogue of the exponential function, since it was suitable to provide a
satisfactory $q$-analogue of the Fourier theory of integral transformations
and series developments. This suitability was made concrete by Bustoz and
Suslov in \cite{BS}, where the authors introduced the subject of $q$-Fourier
series. Some of the subsequent research activity has been already collected
in a book \cite{Suslov}. Among recent developments not yet included in this
book, we quote the orthogonality relations for sums of curly exponential
functions \cite{Koelink}, obtained using spectral methods, and the
construction of a $q$-analogue of the Whittaker-Shannon-Kotel%
\'{}%
nikov sampling theorem \cite{IZa}. The designation "Quantum" has appear
often in recent literature on $q$-analysis, as in the monographs \cite{KC}
and \cite{Ismbook}. This designation is very convenient, since $q$-special
functions are intimately connected with representations of quantum groups 
\cite{QG}.

An abstract formulation designed to capture the essential properties of \ $q$%
-Fourier type systems was proposed in \cite{Ism2001} and we proceed to
describe it here. Let $\mu $ be a measure on the real line and $\{p_{n}(x)\}$
a complete orthonormal system in $L^{2}(\mu )$. Let $(t_{j})$ be a sequence
of points on the real line and assume that $\{r_{n}(x)\}$ is a discrete
orthonormal system with orthogonality relation 
\begin{equation*}
\sum_{j=0}^{\infty }\rho (t_{j})r_{n}(t_{j})\overline{r_{m}(t_{j})}=\delta
_{mn}
\end{equation*}%
and with dual orthogonality 
\begin{equation*}
\sum_{k=0}^{\infty }r_{k}(t_{n})\overline{r_{k}(t_{m})}=\frac{\delta _{mn}}{%
\rho (t_{n})}.
\end{equation*}%
Assume also that the the system $\{r_{n}(x)\}$ is complete in $l^{2}(\sum
\rho (t_{j})\delta _{t_{j}})$. Now define a sequence of functions $%
\{F_{n}(x)\}$ by 
\begin{equation}
F_{n}(x)=\sum_{k=0}^{\infty }r_{k}(t_{n})p_{k}(x)u_{k}  \label{dF}
\end{equation}%
where $\{u_{k}\}$ is an arbitrary sequence of complex numbers in the unit
circle. The following theorem is due to Ismail and comprises in an abstract
form the fundamental fact behind the theory of basic analogs of Fourier
series on a $q$-quadratic grid \cite{BS}.

\textbf{Theorem A \cite{Ism2001} }\emph{The system }$\{F_{n}(x)\}$\emph{\ is
orthogonal and complete in }$L^{2}(\mu )$\emph{.}

To give an idea of what is involved in this statement, we sketch Ismail%
\'{}%
s argument. Since, by (\ref{dF}), $r_{k}(t_{n})$ are the Fourier
coefficients of $F_{n}$ in the basis $\{u_{k}p_{k}\}$, the use of Parseval%
\'{}%
s formula gives 
\begin{equation*}
\int F_{n}(x)\overline{F_{m}(x)}d\mu (x)=\sum_{k=0}^{\infty }r_{k}(t_{n})%
\overline{r_{k}(t_{m})}=\frac{\delta _{mn}}{\rho (x_{n})}
\end{equation*}%
and the orthogonality relation is proved. To show the completeness, choose $%
f\in L^{2}(\mu )$ and assume $\int F_{n}(x)f(x)d\mu (x)=0$ for all $n$.
Again Parseval%
\'{}%
s formula implies $\sum_{k=0}^{\infty }f_{k}\overline{u_{k}r_{k}(t_{m})}=0$
for all $m$, where $f_{k}$ are the Fourier coefficients of $f$ in the basis $%
\{p_{k}\}$. Now the completeness of $\{r_{k}\}$ implies $f_{k}=0$. Therefore 
$f=0$ almost everywhere in $L^{2}(\mu )$.

In \cite{Iprob} (see also section 24.2 of the monograph \cite{Ismbook}),
Ismail posed the problem of studying the mapping properties of operators
with kernels defined as above and conjectured that there was a reproducing
kernel Hilbert space structure behind these operators. We will show that
Ismail%
\'{}%
s conjecture is true. Our approach will reveal a reproducing kernel
structure reminiscent of the well known structure of the Paley-Wiener space
of functions bandlimited to a real interval. However, even in the case when
the system $\{F_{n}(x)\}$ is the set of the complex exponentials, we obtain
results that, as far as our knowledge goes, seem to be new. When the system $%
\{F_{n}(x)\}$ is the set of basis functions of the $q$-Fourier series
constructed with the function $\mathcal{E}_{q}(x;it)$, we will obtain
results that complement the investigations done in \cite{Suslov} and \cite%
{IZa}. In particular it will be shown that a sampling theorem related to the
one derived in \cite{IZa} lives in a reproducing kernel Hilbert space and
that the correspondent $q$-analogues of the $Sinc$ function provide an
orthogonal basis for that space.

Before outlining the paper we want to make clear that the techniques we are
using already exist in some antecedent form. In particular, section two
contains ideas from Ismail theory of generalized $q$-Fourier series and from
Saitoh%
\'{}%
s theory of linear transformations in Hilbert space. However, their
particular combination here leads to new conclusions and sheds new light in
the emerging theory of $q$-Fourier series. It reveals an elementary
structure underlying many systems involving several special functions
simultaneously. The results in section 2.2 and section 3, 4 and 5 were never
explicitly stated before.

The paper can be summarized as follows. The next section contains the
description of the reproducing kernel structure behind the abstract setting
of theorem A. An integral transformation between two Hilbert spaces is
defined in the context of Saitoh%
\'{}%
s theory of linear transformations, basis for both spaces are provided and
the formula for the reproducing kernel of the image Hilbert space is
deduced. In this context a sampling theorem appears, generalizing the one in 
\cite{IZa}. The remaining sections consider three applications of these
results, using specific systems of orthogonal polynomials as well as Bessel
functions and their generalizations. The first application is associated to
formula (\ref{Geg}) and systems of complex exponentials. The second
application is linked to (\ref{qGeg})\ and to systems of curly $q$%
-exponentials, and we write the reproducing kernel as a $_{2}\phi _{2}$
basic hypergeometric function. These two examples explore the interplay
between Lommel polynomials and Bessel functions and the corresponding
relations between their $q$-analogues. In the last section we consider a
construction of a general character, designed originally in the papers \cite%
{IZ},\ \cite{IRZ} and \cite{Ism2001}. It allows to extend the interplay
between Bessel functions and Lommel polynomials to a more general class of
functions. Using this construction we will make a brief discussion about the
application of our results to spaces weighted by Jacobi weights and their $q$%
-analogues.

\section{The reproducing kernel structure}

Let $H$ be a class of complex valued functions, defined in a set $X\subset 
\mathbf{C}$, such that $H$ is a Hilbert space with the norm of $L^{2}\left(
X,\mu \right) $. The function $k\left( s,x\right) $ is a \emph{reproducing
kernel} of $H$ if

$i)$ $k\left( .,x\right) \in H$ for every $x\in X;$

$ii)$ $f\left( x\right) =\left\langle f\left( .\right) ,k\left( .,x\right)
\right\rangle $ for every $f\in H$, $x\in X.$

The space $H$ is said to be a reproducing kernel Hilbert space if it
contains a reproducing kernel. From a structural point of view, the correct
approach to the study of our problem is via Saitoh%
\'{}%
s theory of linear transforms of Hilbert space.

\subsection{Preliminaries on Saitoh%
\'{}%
s theory of linear transformations in Hilbert space}

This theory can be found in works by Saitoh \cite{S1}, \cite{S2} and we
proceed to give a brief account of the results that we are going to use. An
account of the results quoted in this section can also be found in Higgins
recent survey \cite{H2}.

For each $t$ belonging to a domain $D$, let $K_{t}$ belong to $H$ (a
separable Hilbert space). Then, 
\begin{equation*}
k(t,s)=\left\langle K_{t},K_{s}\right\rangle _{H}
\end{equation*}%
is defined on $D\times D$ and is called the kernel function of the map $%
K_{t} $. Now consider the set of images of $H$ by the transformation 
\begin{equation*}
(Kg)(t)=\left\langle g,K_{t}\right\rangle _{H}=f(t)
\end{equation*}%
and denote this set of images by $R(K)$. The following theorem can be found
in \cite{S1}:

\textbf{Theorem B }\emph{The kernel }$k(t,s)$\emph{\ determines uniquely a
reproducing kernel Hilbert space for which it is the reproducing kernel.
This reproducing kernel Hilbert space is precisely }$R(K)$\emph{\ and it can
have no other reproducing kernel.}

Now, suppose that $\{K_{t}\}(t\in D)$ is complete in $H$ so that $K$ is one
to one. Then we have 
\begin{equation*}
\left\Vert Kg\right\Vert _{R(K)}=\left\Vert g\right\Vert _{H}.
\end{equation*}%
The following theorem, due to Higgins \cite{Higgins}, will be critical on
the remainder. We will use only the following "orthogonal basis case", a
special case of the theorem in \cite{Higgins}:

\textbf{Theorem C }\emph{With the notations established earlier, we have: If
there exists }$\left\{ t_{n}\right\} (n\in \mathbf{I\subset Z})$ \emph{such
that }$\left\{ K_{t_{n}}\right\} $\emph{\ is an orthogonal basis, we then
have the sampling expansion } 
\begin{equation*}
f(t)=\sum_{n\in \mathbf{I}}f(t_{n})\frac{k(t,t_{n})}{k(t_{n},t_{n})}
\end{equation*}%
\emph{in }$R(K)$\emph{, pointwise over }$\mathbf{I}$\emph{, and uniformly
over any compact subset of }$D$\emph{\ for which }$\left\Vert
K_{t}\right\Vert $\emph{\ is bounded.}

\subsection{The reproducing kernel for $q$-Fourier type systems}

In this section we will show the existence of a reproducing kernel structure
behind the abstract setting of Theorem A. The results will follow from the
study of the mapping properties of an integral transform whose kernel is
obtained from the sequence of functions $\{r_{k}\}$ and $\{p_{k}\}$.

\begin{theorem}
There exists a kernel $K(x,t)$ satisfying the requirements of Saitoh%
\'{}%
s theory of linear transformations, such that $K(x,x_{n})=\lambda
_{n}F_{n}(x)$, where $\{F_{n}(x)\}$ is the orthogonal sequence of functions
in Theorem A and $\{\lambda _{n}\}$ is a sequence of real numbers.
\end{theorem}

\begin{proof}
Our first technical problem comes from the fact that, when $\{r_{k}\}$ is a
discrete system of orthogonal polynomials with a determinate moment problem,
then $\{r_{k}(t)\}\in l^{2}$ if and only if $x$ is a mass point for the
measure of orthogonality. For this reason the series 
\begin{equation*}
\sum_{k=0}^{\infty }r_{k}(t)p_{k}(x)u_{k}
\end{equation*}%
would diverge if $t$ is not such a point (this is pointed out in Section 5
of \cite{Ism2001}). Since we want our kernel to be defined for every $t$, we
will assume the existence of an auxiliary system of functions $\{\mathcal{J}%
_{k}(x)\}\in l^{2}$ for every $t$ real, and such that every function $%
\mathcal{J}_{k}$ interpolates $r_{k}$ at the mass points $\{x_{n}\}$ in the
sense that 
\begin{equation}
\mathcal{J}_{k}(x_{n})=\lambda _{n}r_{k}(\frac{1}{x_{n}})  \label{intgeral}
\end{equation}%
for every $k=0,1,...$ and $n=0,1,..$.and some constant $\lambda _{n}$
independent of $k$.\ Now we can use the functions $\mathcal{J}_{k}(t)$ to
define a kernel $K(x,t)$ as 
\begin{equation}
K(x,t)=\sum_{k=0}^{\infty }\mathcal{J}_{k}(t)p_{k}(x)u_{k}.  \label{defk}
\end{equation}%
Such a kernel is well defined and belongs to $L^{2}(\mu )$, since it is a
sum of basis functions of $L^{2}(\mu )$. From (\ref{dF}), (\ref{intgeral})\
and (\ref{defk}) we have 
\begin{eqnarray*}
K(x,x_{n}) &=&\sum_{k=0}^{\infty }\mathcal{J}_{k}(x_{n})p_{k}(x)u_{k} \\
&=&\lambda _{n}\sum_{k=0}^{\infty }r_{k}(\frac{1}{x_{n}})p_{k}(x)u_{k} \\
&=&\lambda _{n}F_{n}(x).
\end{eqnarray*}
\end{proof}

\begin{remark}
Observe that theorem 1 and theorem A with $t_{n}=\frac{1}{x_{n}}$ show that $%
K(x,x_{n})$ is an orthogonal basis for the space $L^{2}(\mu )$.
\end{remark}

\begin{remark}
In the abstract formulation it may not be clear why the constant $\lambda
_{n}$ must be present. Actually the construction would work without it, but
for technical reasons that will become evident upon consideration of
examples we prefer to use it. Otherwise, careful bookkeeping of the
normalization constants would be required in the remaining sections.
\end{remark}

\begin{remark}
It will be seen in the last section that a general constructive method is
available in order to find the function $\mathcal{J}_{k}$ under very natural
requirements on the polynomials $r_{k}$.
\end{remark}

Now define an integral transformation $F$ by setting 
\begin{equation*}
(Ff)(t)=\int f(x)K(x,t)d\mu (x).
\end{equation*}%
We will study this transform as a map whose domain is the Hilbert space $%
L^{2}(\mu )$. Endowing the range of $F$ with the inner product 
\begin{equation}
\left\langle Ff,Fg\right\rangle _{F(L^{2}(\mu ))}=\left\langle
f,g\right\rangle _{L^{2}(\mu )}  \label{inn}
\end{equation}%
then $F(L^{2}(\mu ))$ becomes a Hilbert space isometrically isomorphic to $%
L^{2}(\mu )$ under the isomorphism $F$. Using Saitoh%
\'{}%
s theory with $D=\mathbf{R}$, $K_{t}=K(.,t)$, $H=L^{2}(\mu )$ and $%
R(K)=F(L^{2}(\mu ))$, we obtain at once:

\begin{theorem}
The transform $F$ is a Hilbert space isomorphism\ mapping the space $%
L^{2}(\mu )$ into $F(L^{2}(\mu ))$. The space $F(L^{2}(\mu ))$ is a Hilbert
space with reproducing kernel given by 
\begin{equation}
k(t,s)=\int K(x,t)\overline{K(x,s)}d\mu (x).  \label{rep}
\end{equation}
\end{theorem}

An interesting feature of this particular setting is that every function in $%
F(L^{2}(\mu ))$ has two different expansions: One is the sampling expansion
naturally associated with the reproducing kernel structure, the other is the
expansion in the basis $\mathcal{J}_{n}(x)$. It should be remarked that, in
most of the previously known sampling expansions, these two expansions were
the same. These expansion results are summarized in the next theorem.

\begin{theorem}
Every function $f$ of the form 
\begin{equation}
f(x)=\int u(t)K(t,x)d\mu (t)  \label{absbandlim}
\end{equation}%
with $u\in L^{2}(\mu )$, admits an expansion 
\begin{equation}
f(t)=\sum_{n=0}^{\infty }a_{n}\mathcal{J}_{n}(t)  \label{expansion}
\end{equation}%
where the coefficients $a_{k}$ are given by 
\begin{equation*}
a_{n}=u_{n}\left\langle u,p_{n}(.)\right\rangle _{L^{2}(\mu )}
\end{equation*}%
and a sampling expansion 
\begin{equation}
f(x)=\sum f(t_{n})\frac{k(x,t_{n})}{k(t_{n},t_{n})}.  \label{samp}
\end{equation}%
The sum in (\ref{samp}) converges absolutely. Furthermore, it converges
uniformly in every set such that $\left\Vert K(.,t)\right\Vert _{L^{2}(\mu )}
$ is finite.
\end{theorem}

\begin{proof}
We already know by default that $\{p_{n}(x)\}$ is a basis for $L^{2}(\mu )$.
It remains to prove that $\{\mathcal{J}_{n}(t)\}$ is a basis for $%
F(L^{2}(\mu ))$. Observe that 
\begin{eqnarray*}
(Fp_{n})(t) &=&\int p_{n}(x)K(x,t)d\mu (x) \\
&=&\sum_{k=0}^{\infty }\mathcal{J}_{k}(t)u_{k}\int p_{n}(x)p_{k}(x)d\mu (x)
\\
&=&\mathcal{J}_{n}(t)u_{n}.
\end{eqnarray*}%
Since $\{p_{n}(x)\}$ is a basis for $L^{2}(\mu )$ and $F$ is an isomorphism
between $L^{2}(\mu )$ and $F(L^{2}(\mu ))$, then $\{u_{n}\mathcal{J}_{n}(x)\}
$ is a basis for $F(L^{2}(\mu ))$. To prove the last assertion of the
theorem, observe that the function $f$ defined by (\ref{absbandlim}) belongs
to $F(L^{2}(\mu ))$ and therefore can be expanded in the basis $\{u_{n}%
\mathcal{J}_{n}(x)\}$. The Fourier coefficients of this expansion are 
\begin{equation*}
a_{n}=\left\langle f,u_{n}\mathcal{J}_{n}(.)\right\rangle _{F(L^{2}(\mu
))}=\left\langle Fu,F(p_{n}(.))\right\rangle _{F(L^{2}(\mu ))}=\left\langle
u,p_{n}(.)\right\rangle _{L^{2}(\mu )}
\end{equation*}%
where we have used (\ref{inn}) in the last identity. The sampling expansion
follows from applying theorem C to our setting and using remark 1.
\end{proof}

\begin{remark}
The construction of this section has never appeared before in the
literature, but it is reminiscent of the reproducing kernel structure of the
Paley-Wiener space. In the classical situations (see for example \cite%
{Garcia} and \cite{NW}\ for an account of these constructions with several
examples) generalizing this structure, there is an integral transform whose
kernel is defined as 
\begin{equation}
K(x,t)=\sum S_{k}(t)e_{k}(x)  \label{kerclass}
\end{equation}%
where $e_{k}(x)$ is an orthogonal basis for the domain Hilbert space and $%
S_{k}(t)$ is a sequence of functions such that there exists a sequence $%
\{t_{n}\}$ satisfying the sampling property 
\begin{equation}
S_{k}(t_{n})=a_{n}\delta _{n,k}  \label{sp}
\end{equation}%
As an instance, take $S_{k}(t)=\frac{\sin \pi (t-k)}{\pi (t-k)}$ and $%
e_{k}(x)=e^{ikx}$. Then (\ref{kerclass}) is 
\begin{equation*}
e^{itx}=\sum_{k=-\infty }^{\infty }\frac{\sin \pi (t-k)}{\pi (t-k)}e^{ikx}
\end{equation*}%
and $K(x,t)$ is the kernel of the Fourier transform. The corresponding
reproducing kernel Hilbert space is the Paley-Wiener space. The root of
these ideas is in Hardy%
\'{}%
s groundbreaking paper \cite{Hardy}. For an application of this classical
set up to Jackson $q$-integral transforms and the third Jackson $q$-Bessel
function, see \cite{Abr}. In our construction we made a modification of this
classical setting: Instead of the sequence of functions $S_{k}$, with the
sampling property (\ref{sp}), we considered a sequence of functions $\{%
\mathcal{J}_{k}\}$, interpolating an orthogonal system $\{r_{k}\}$ in the
sense of\ (\ref{intgeral}). And we have seen that the essential properties
of classical reproducing kernel settings are kept. However, this
modification allows to recognize a class of reproducing kernel Hilbert
spaces that were obscured until now. This will become clear in the next
section.
\end{remark}

\section{The Fourier system with ultraspherical weights}

The $nth$ ultraspherical (or Gegenbauer) polynomial of order $\nu $ is
denoted by $C_{n}^{\nu }(x)$. These polynomials satisfy the orthogonality
relation 
\begin{equation*}
\int_{-1}^{1}C_{n}^{\nu }(x)C_{m}^{\nu }(x)(1-x^{2})^{\nu -1/2}dx=\frac{%
(2\nu )_{n}\sqrt{\pi }\Gamma (\nu +\frac{1}{2})}{n!(\nu +n)\Gamma (\nu )}%
\delta _{m,n}
\end{equation*}%
and form a complete sequence in the Hilbert space $L^{2}[(-1,1),(1-x^{2})^{%
\nu -1/2}]$. For typographical convenience we will introduce the following
notation for this Hilbert space: 
\begin{equation*}
H^{\nu }=L^{2}[(-1,1),(1-x^{2})^{\nu -1/2}].
\end{equation*}%
The Bessel function of order $\nu $, $J_{\nu }(x)$, is defined by the power
series expansion 
\begin{equation}
J_{\nu }(z)=\sum\limits_{n=0}^{\infty }\frac{(-1)^{n}}{n!\Gamma (\nu +n+1)}%
\left( \frac{z}{2}\right) ^{\nu +2n}.  \label{Bessel}
\end{equation}%
The $nth$ Lommel polynomial of order $\nu $, denoted by $h_{n,\nu }(x)$, is
related to the Bessel functions by the relation 
\begin{equation}
J_{\nu +k}(x)=h_{k,\nu }(\frac{1}{x})J_{\nu }(x)-h_{k-1,\nu -1}(\frac{1}{x}%
)J_{\nu -1}(x)\text{.}  \label{besslom}
\end{equation}%
The Lommel polynomials satisfy the discrete orthogonality relation 
\begin{equation*}
\sum_{k=0}^{\infty }\frac{1}{(j_{\nu ,k})^{2}}h_{n,\nu +1}(\pm \frac{1}{%
j_{\nu ,k}})h_{m,\nu +1}(\pm \frac{1}{j_{\nu ,k}})=\frac{1}{2(\nu +n+1)}%
\delta _{nm}
\end{equation*}%
and the dual orthogonality 
\begin{equation*}
\sum_{k=0}^{\infty }2(\nu +n+1)h_{k,\nu +1}(\pm \frac{1}{j_{\nu ,n}}%
)h_{k,\nu +1}(\pm \frac{1}{j_{\nu ,m}})=(j_{\nu ,k})^{2}\delta _{nm}.
\end{equation*}%
They form a complete orthogonal system in the $l^{2}$ space weighted by the
discrete measure with respect to which they are orthogonal. We will use
these two complete orthogonal systems in our first illustration of the
general results. Set 
\begin{equation*}
p_{k}(x)=\sqrt{\frac{k!(\nu +k)}{(2\nu )_{k}}}C_{k}^{\nu }(x)
\end{equation*}%
and 
\begin{equation*}
r_{k}(t)=\sqrt{2(\nu +n)}h_{k,\nu -1}(t).
\end{equation*}%
Consider also 
\begin{equation*}
\mathcal{J}_{k}(t)=\sqrt{2(\nu +n)}J_{\nu +k}(t).
\end{equation*}%
Denote by $j_{\nu ,k}$\ the $kth$ zero of the Bessel function of order $\nu $%
. Substituting $x=j_{\nu ,n}$ in (\ref{besslom}), the following
interpolating property is obtained 
\begin{equation}
h_{k,\nu -1}(\frac{1}{j_{\nu ,n}})=-\frac{J_{\nu +k}(j_{\nu ,n})}{J_{\nu
-1}(j_{\nu ,n})}.  \label{int}
\end{equation}%
The interpolating property (\ref{int}) will play the role of (\ref{intgeral}%
) with $\lambda _{n}=-\frac{1}{J_{\nu -1}(j_{\nu ,n})}$. Consider also the
sequence of complex numbers $\{u_{n}\}$ defined as 
\begin{equation*}
u_{k}=i^{k}
\end{equation*}%
and set 
\begin{equation*}
K^{\nu }(x,t)=\sqrt{2}\sum_{k=0}^{\infty }i^{k}(\nu +k)\sqrt{\frac{2k!}{%
(2\nu )_{k}}}J_{\nu +k}(t)C_{k}^{\nu }(x).
\end{equation*}%
Now, theorem 1 tells us that $K^{\nu }(x,j_{\nu ,n})$ is an orthogonal basis
of the space $H^{\nu }$. Moreover, formula (\ref{Geg}) implies%
\begin{equation*}
\sqrt{\frac{\pi t}{2}}K^{\frac{1}{2}}(x,t)=e^{ixt}
\end{equation*}%
and therefore we can think of $\Gamma (\nu )\left( \frac{t}{2}\right) ^{\nu
}K^{\nu }(x,t)$ as a one parameter generalization of the complex exponential
kernel that may be worth of further study. In the special case $\nu =\frac{1%
}{2}$ we recover the well known orthogonality and completeness of the
complex exponentials $\{e^{i\pi nx}\}$ in $L^{2}(-1,1)$.

The transformation $F$ is defined, for every $f\in H^{\nu }$, as 
\begin{equation*}
(Ff)(t)=\int_{-1}^{1}f(x)K^{\nu }(x,t)(1-x^{2})^{\nu -1/2}dx
\end{equation*}%
and theorem 2 gives that the reproducing kernel of $F(H^{\nu })$ is 
\begin{equation*}
k^{\nu }(t,s)=\int_{-1}^{1}K^{\nu }(x,t)\overline{K^{\nu }(x,s)}%
(1-x^{2})^{\nu -1/2}dx.
\end{equation*}%
When $\nu =1/2$ this becomes 
\begin{equation*}
k^{\frac{1}{2}}(t,s)=\frac{\pi \sqrt{ts}}{2}\frac{\sin (t-s)}{(t-s)}.
\end{equation*}%
Since $\{\sqrt{\frac{k!(\nu +k)}{(2\nu )_{k}}}C_{n}^{\nu }(t)\}$ forms a
basis of the space $H^{\nu }$, then also $\{i^{n}\sqrt{2(\nu +n)}J_{\nu
+n}(x)\}=F\{\sqrt{\frac{k!(\nu +k)}{(2\nu )_{k}}}C_{n}^{\nu }(t)\}$ is a
basis of the space $F\left( H^{\nu }\right) $. In this situation theorem 3
reads:

\begin{theorem}
Let $f$ be a function of the form 
\begin{equation}
f(t)=\int_{-1}^{1}u(x)K^{\nu }(x,t)(1-x^{2})^{\nu -1/2}dx  \label{utrbandlim}
\end{equation}%
where $u\in H^{\nu }$. Then $f$ can be written as 
\begin{equation}
f(t)=\sum_{n=0}^{\infty }a_{n}J_{\nu +n}(t)  \label{Neumann}
\end{equation}%
with the coefficients $a_{n}$ given by 
\begin{equation}
a_{n}=i^{n}(\nu +n)\sqrt{\frac{2n!}{(2\nu )_{n}}}\int_{-1}^{1}u(x)C_{n}^{\nu
}(x)(1-x^{2})^{\nu -1/2}dx  \label{coef}
\end{equation}%
and also as the sampling formula%
\begin{equation*}
f(t)=\sum_{n=0}^{\infty }f(j_{\nu ,n})\frac{k^{\nu }(t,j_{\nu ,n})}{k^{\nu
}(j_{\nu ,n},j_{\nu ,n})}.
\end{equation*}
\end{theorem}

\begin{remark}
Expansions of the type (\ref{Neumann}) are known as Neumann series of Bessel
functions (see chapter 16 of \cite{Watson}).
\end{remark}

\begin{remark}
When $\nu =1/2$ the above sampling theorem states that every function of the
form 
\begin{equation*}
f(t)=\left( \frac{t}{2}\right) ^{-\frac{1}{2}}\int_{-1}^{1}u(x)e^{ixt}dx,
\end{equation*}%
with $u\in L^{2}[(-1,1)],$ can be represented as 
\begin{equation*}
f(t)=\sum_{n=0}^{\infty }f(2\pi n)\sqrt{\frac{t}{2\pi n}}\frac{\sin (t-2\pi
n)}{(t-2\pi n)}.
\end{equation*}
\end{remark}

\section{The $q$-Fourier system with $q$-ultraspherical weights}

We proceed to describe the $q$-analogue of the previous situation. Choose a
number $q$ such that $0<q<1$. The now classical notational conventions from 
\cite{GR} and \cite{Ismbook} for $q$-infinite products and basic
hypergeometric series will be used often.

The $q$-exponential function that we talked about in the introduction is
defined in terms of basic hypergeometric series as 
\begin{equation*}
\mathcal{E}_{q}(x;t)=\frac{(-t;q^{\frac{1}{2}})_{\infty }}{%
(qt^{2};q^{2})_{\infty }}{}{}_{2}\phi _{1}\left( \left. 
\begin{array}{c}
q^{\frac{1}{4}}e^{i\theta },q^{\frac{1}{4}}e^{-i\theta } \\ 
q^{\frac{1}{2}}%
\end{array}%
\right\vert q^{\frac{1}{2}},-t\right) 
\end{equation*}%
where $x=\cos \theta $. The continuous $q$-ultraspherical polynomials of
order $\nu $ are denoted by $C_{n}^{\nu }(x;q^{\nu }|q)$ and satisfy the
orthogonality 
\begin{equation*}
\int_{-1}^{1}C_{n}^{\nu }(x;q^{\nu }|q)C_{m}^{\nu }(x;q^{\nu }|q)w(x;q^{\nu
}\mid q)dx=\frac{(1-q^{\nu })(q^{2\nu };q)_{n}}{(1-q^{n+\nu })(q;q)_{n}}%
\delta _{m,n}
\end{equation*}%
where the weight function $w(x;\beta \mid q)$ is 
\begin{equation*}
w(\cos \theta ;\beta |q)=\frac{(q,q^{2\nu };q)_{\infty }(e^{2i\theta
},e^{-2i\theta };q)_{\infty }}{\sin \theta (2\pi q^{\nu },q^{\nu
+1};q)_{\infty }(\beta e^{2i\theta },\beta e^{-2i\theta };q))_{\infty }}%
,(0<\theta <\pi ).
\end{equation*}%
The polynomials $\{C_{n}^{\nu }(x;q^{\nu }|q)\}$\ form a basis of the
Hilbert space $H_{q}^{\nu }$ defined as 
\begin{equation*}
H_{q}^{\nu }=L^{2}[(-1,1),w(x;q^{\nu }\mid q)].
\end{equation*}%
The second Jackson $q$-Bessel function of order $\nu $ is defined by the
power series 
\begin{equation*}
J_{\nu }^{(2)}(x;q)=\frac{(q^{\nu +1};q)_{\infty }}{(q;q)_{\infty }}%
\sum_{k=0}^{\infty }(-1)^{n}\frac{(x/2)^{\nu +2n}}{(q;q)_{n}(q^{\nu
+1};q)_{n}}q^{n(\nu +n)}.
\end{equation*}%
Since this is the only $q$-Bessel function to be used in the text, we will
drop the superscript for shortness of the notation and write $J_{\nu
}(x;q)=J_{\nu }^{(2)}(x;q)$. The $q$-Lommel polynomials associated to the
Jackson $q$-Bessel function of order $\nu $ are denoted by $h_{n,\nu -1}(x;q)
$. These polynomials were defined in \cite{Ism82} by means of the relation 
\begin{equation}
q^{n\nu +n(n-1)/2}J_{\nu +n}(x;q)=h_{n,\nu }(\frac{1}{x};q)J_{\nu
}(x;q)-h_{n-1,\nu -1}(\frac{1}{x};q)J_{\nu -1}(x;q).  \label{qbesslom}
\end{equation}%
The $q$-Lommel polynomials satisfy the orthogonality relation, for a certain
constant $A_{n}(\nu +1)$ which is not explicitly known (see section 14.4 of 
\cite{Ismbook}):%
\begin{equation*}
\sum_{n=1}^{\infty }\frac{A_{n}(\nu +1)}{(j_{\nu ,n}(q))^{2}}h_{n,\nu
+1}(\pm \frac{1}{j_{\nu ,n}(q)};q)h_{m,\nu +1}(\pm \frac{1}{j_{\nu ,n}(q)}%
;q)=\frac{q^{n\nu +n(n+1)/2}}{1-q^{n+\nu +1}}\delta _{nm}
\end{equation*}%
and the dual orthogonality 
\begin{equation*}
\sum_{n=1}^{\infty }\frac{(1-q^{n+\nu +1})}{q^{n\nu +n(n+1)/2}}h_{n,\nu
+1}(\pm \frac{1}{j_{\nu ,n}(q)};q)h_{m,\nu +1}(\pm \frac{1}{j_{\nu ,n}(q)}%
;q)=\frac{(j_{\nu ,n}(q))^{2}}{A_{k}(\nu +1)}\delta _{nm}.
\end{equation*}%
Consider 
\begin{equation*}
p_{k}(x)=\sqrt{\frac{(1-q^{k+\nu })(q;q)_{k}}{(1-q^{\nu })(q^{2\nu };q)_{k}}}%
C_{k}(x;q^{\nu }|q),
\end{equation*}%
\begin{equation*}
r_{k}(t)=\sqrt{\frac{(1-q^{k+\nu })}{(1-q^{\nu })}}q^{-k\nu
/2-k(k-1)/4}h_{k,\nu }(2t;q)
\end{equation*}%
and 
\begin{equation*}
\mathcal{J}_{k}(t)=\sqrt{\frac{(1-q^{k+\nu })}{(1-q^{\nu })}}q^{k\nu /2+%
\frac{k(k-1)}{4}}J_{\nu +k}(2t;q).
\end{equation*}%
The parameters $u_{k}$ will be given by 
\begin{equation*}
u_{k}=i^{k}.
\end{equation*}%
Denote by $j_{\nu ,k}(q)$ the $kth$ zero of $J_{\nu }(x;q)$.\ Setting $%
t=j_{\nu ,k}(q)$ in (\ref{qbesslom}) we have the interpolating property 
\begin{equation}
h_{k,\nu -1}(\frac{1}{j_{\nu ,n}(q)};q)=-q^{k\nu +k(k-1)/2}\frac{J_{\nu
+k}(j_{\nu ,n}(q);q)}{J_{\nu -1}(j_{\nu ,n}(q);q)}.  \label{qint}
\end{equation}%
This means that in (\ref{intgeral}) we must take $\lambda _{n}=-\frac{1}{%
J_{\nu -1}(j_{\nu ,n}(q);q)}$. In this context, the kernel $K(x,t)$ is given
as 
\begin{equation*}
K_{q}^{\nu }(x,t)=\sum_{k=0}^{\infty }i^{k}q^{k(2\nu +k-1)/4}\frac{%
(1-q^{k+\nu })}{(1-q^{\nu })}\sqrt{\frac{(q;q)_{k}}{(q^{2\nu };q)_{k}}}%
J_{\nu +k}(2t;q)C_{k}(x;q^{\nu }|q)
\end{equation*}%
and the use of (\ref{qGeg}) gives, when $\nu =\frac{1}{2}$,%
\begin{equation}
K_{q}^{\frac{1}{2}}(x,t)=\frac{(-qt^{2};q^{2})_{\infty }(q^{\nu
+1};q)_{\infty }}{(q;q)_{\infty }}t^{\nu }\mathcal{E}_{q}(x;it).
\label{qexp}
\end{equation}%
The basis functions of the domain space are 
\begin{equation*}
F_{n}(x)=K_{q}^{\nu }(x,j_{\nu ,k}(q)).
\end{equation*}%
When $\nu =\frac{1}{2}$ this gives that $\{\mathcal{E}_{q}(x;ij_{\frac{1}{2}%
,n-1}(q))\}$ is orthogonal and complete in $H_{q}^{\nu }$. This is the case
of the $q$-Fourier series studied in \cite{Suslov}. Now define the transform 
\begin{equation}
(F_{q}^{\nu }f)(t)=\int_{-1}^{1}f(x)K_{q}^{\nu }(x,t)w(x;q^{\nu }\mid q)dx
\label{qbl}
\end{equation}%
for every $f\in H_{q}^{\nu }$. We have 
\begin{equation*}
i^{k}\sqrt{\frac{(1-q^{k+\nu })}{(1-q^{\nu })}}q^{k\nu /2+\frac{k(k-1)}{4}%
}J_{\nu +k}(2t;q)=F_{q}^{\nu }\left( \sqrt{\frac{(1-q^{k+\nu })(q;q)_{k}}{%
(1-q^{\nu })(q^{2\nu };q)_{k}}}C_{k}(x;q^{\nu }|q)\right) 
\end{equation*}%
and $\{J_{\nu +n}(t;q)\}_{n=0}^{\infty }$ is a basis of the space $%
F_{q}^{\nu }\left( H_{q}^{\nu }\right) $. We can also state the $q$-analogue
of the expansion in theorem 5 and obtain a $q$-Neumann expansion theorem in $%
q$-Bessel functions.

\begin{theorem}
Let $f$ be a function of the form 
\begin{equation}
f(t)=\int_{-1}^{1}u(x)K_{q}^{\nu }(x,t)w(x;q^{\nu }\mid q)dx
\label{qutrbandlim}
\end{equation}%
where $u\in H_{q}^{\nu }$. Then $f$ can be written as 
\begin{equation*}
f(t)=\sum_{n=0}^{\infty }a_{n}J_{\nu +n}(t;q)
\end{equation*}%
with the coefficients $a_{n}$ given by 
\begin{equation*}
a_{n}=i^{n}\frac{(1-q^{k+\nu })}{(1-q^{\nu })}\sqrt{\frac{(q;q)_{k}}{%
(q^{2\nu };q)_{k}}}\int_{-1}^{1}u(x)C_{n}(x;q^{\nu }|q)w(x;q^{\nu }\mid q)dx.
\end{equation*}
\end{theorem}

Once more we know by theorem 2 that the space $F_{q}^{\nu }\left( H_{q}^{\nu
}\right) $ is a space with reproducing kernel $k_{q}^{\nu }(t,s)$, given by%
\begin{equation*}
k_{q}^{\nu }(t,s)=\int_{-1}^{1}K_{q}^{\nu }(x,t)\overline{K_{q}^{\nu }(x,s)}%
w(x;q^{\nu }\mid q)dx.
\end{equation*}%
When $\nu =\frac{1}{2}$ we can compute the reproducing kernel in an explicit
form using the following integral, evaluated in \cite{IS}: 
\begin{eqnarray}
&&\int_{0}^{\pi }\mathcal{E}_{q}(\cos \theta ;\alpha )\mathcal{E}_{q}(\cos
\theta ;\beta )\frac{(e^{2i\theta },e^{-2i\theta };q)_{\infty }}{(\gamma
e^{2i\theta },\gamma e^{-2i\theta };q)_{\infty }}\;d\theta  \label{Isf} \\
&{}&\quad =\frac{2\pi (\gamma ,q\gamma ,-\alpha \beta q^{1/2};q)_{\infty }}{%
(q,\gamma ^{2};q)_{\infty }(q\alpha ^{2},q\beta ^{2};q^{2})_{\infty }}{}%
_{2}\phi _{2}\left( \left. 
\begin{array}{c}
-q^{1/2}\alpha /\beta ,-q^{1/2}\beta /\alpha \\ 
q\gamma ,,-\alpha \beta \gamma q^{1/2}%
\end{array}%
\right\vert q,-\alpha \beta \gamma q^{1/2}\right) .  \notag
\end{eqnarray}

\begin{theorem}
The space $F_{q}^{\frac{1}{2}}\left( H_{q}^{\frac{1}{2}}\right) $ is a space
with reproducing kernel $k_{q}^{\frac{1}{2}}(t,s)$, given by%
\begin{equation*}
k_{q}^{\nu }(t,s)=2\pi \left[ \frac{(q^{3/2};q)_{\infty }^{3/2}}{%
(q;q)_{\infty }^{2}}\right] ^{2}
\end{equation*}
\begin{equation}
\times (q^{1/2},-tsq^{1/2};q)_{\infty }(ts)^{\frac{1}{2}}{}_{2}\phi
_{2}\left( \left. 
\begin{array}{c}
-q^{1/2}t/s,-q^{1/2}s/t \\ 
q^{3/2},-tsq%
\end{array}%
\right\vert q,-tsq\right) .  \label{qrp}
\end{equation}
\end{theorem}

\begin{proof}
Applying theorem 2 and (\ref{qexp})\ we know that $k_{q}^{\nu }(t,s)$ is
given by 
\begin{equation*}
k_{q}^{\nu }(t,s)=\left[ \frac{(q^{\nu +1};q)_{\infty }}{(q;q)_{\infty }}%
\right] ^{2}(-qt^{2},-qs^{2};q^{2})_{\infty }(ts)^{\frac{1}{2}}
\end{equation*}%
\begin{equation*}
\times \int_{-1}^{1}\mathcal{E}_{q}(x;it)\mathcal{E}_{q}(x;-is)w(x;q^{\nu
}\mid q)dx.
\end{equation*}%
Make the substitutions $x=\cos \theta $, $it=\alpha $, $is=\beta $, and $%
q^{\nu }=\gamma $ in (\ref{Isf}). Then (\ref{qrp}) follows.
\end{proof}

As in the preceding sections we can formulate a sampling theorem.

\begin{theorem}
Every function of the form (\ref{qutrbandlim}) admits the expansion 
\begin{equation}
f(x)=\sum_{k=0}^{\infty }f(t_{k})\frac{k_{q}^{\nu }(x,t_{k})}{k_{q}^{\nu
}(t_{k},t_{k})}  \label{qsamp}
\end{equation}%
where $t_{k}=\frac{j_{\nu ,k}(q)}{2}.$
\end{theorem}

\begin{remark}
The functions $\frac{k_{q}^{\nu }(x,t_{k})}{k_{q}^{\nu }(t_{k},t_{k})}$ are
related to the functions $\func{Si}nc_{q}(t,k)$ from \cite{IZa}. In \cite%
{IZa} it is shown that there is a sampling theorem associated to the $\func{%
Si}nc_{q}(t,k)$ functions that can be written as an interpolating formula of
the Lagrange type. The above discussion shows what is the sampling theorem
associated to the reproducing kernel Hilbert space and its relation to the
one obtained in \cite{IZa}.
\end{remark}

\begin{remark}
By the proof of Theorem B (see \cite{H2}), the functions $\frac{k_{q}^{\nu
}(x,t_{k})}{k_{q}^{\nu }(t_{k},t_{k})}$are orthogonal in the space $F_{q}^{%
\frac{1}{2}}(H_{q}^{\frac{1}{2}})$. This result is a $q$-analogue of the
important fact, proved by Hardy in \cite{Hardy}, that the functions $\frac{%
\sin \pi (t-k)}{\pi (t-k)}$ are orthogonal in the classical Paley-Wiener
space.
\end{remark}

\begin{remark}
Important information concerning the zeros of the second Jackson $q$-Bessel
function, that appear as sampling nodes in the expansion (\ref{qsamp}), was
obtained very recently by Walter Hayman in \cite{Hayman} using a method due
to Bergweiler and Haymann \cite{BH}: the asymptotic expansion 
\begin{equation*}
j_{\nu ,k}^{2}(q)=4q^{1-2n-\nu }\{1+\sum_{\nu =1}^{n}b_{\nu }q^{k\nu
}+O\left| q^{(n+1)k}\right| \}
\end{equation*}
as $k\rightarrow \infty $,\ with the constants $b_{\nu }$ depending on $a$
and $q$, holds. Therefore, for big $k$, the sampling nodes are exponentially
separated in a similar way to what was observed in \cite{Abr} and \cite{ABC}%
. In the case where $\nu =\pm \frac{1}{2}$, the zeros were studied by Suslov 
\cite{Sus}.
\end{remark}

\section{A generalization}

We begin this last section describing a formal approach generalizing the
situations studied in the two previous sections. This formal approach was
initiated in \cite{IZ} and \cite{IRZ} with the purpose of finding functions
to play the role of the Lommel polynomials in more general situations, and
was studied further in \cite{Ism2001}. In the context studied in this paper,
it will be of particular relevance, since it gives a constructive method to
find the functions $\mathcal{J}_{k}$ satisfying (\ref{intgeral}). Let $%
\{f_{n,\nu }\}$ be a sequence of polynomials defined recursively by $%
f_{0,\nu }(x)=1$, $f_{1,\nu }(x)=xB_{\nu }$ and 
\begin{equation*}
f_{n+1,\nu }(x)=[xB_{n+\nu }]f_{n,\nu }(x)-C_{n+\nu -1}f_{n-1,\nu }(x).
\end{equation*}%
Assuming the positivity condition $B_{n+\nu }B_{n+\nu +1}C_{n+\nu }>0$ $%
(n\geq 0)$ and the convergence of the series $\sum_{n=0}^{\infty }\frac{%
C_{n+\nu }}{B_{n+\nu }B_{n+\nu +1}}$, it can be shown, using facts from the
general theory of orthogonal polynomials, that the polynomials $f_{n,\nu }$
are orthogonal with respect to a compact supported discrete measure and that
the support points of this measure are $\frac{1}{x_{n,\nu }}$, where the $%
x_{n,\nu }$ are the zeros of an entire function $\mathcal{J}$ satisfying 
\begin{equation}
C_{\nu }...C_{\nu +n-1}\mathcal{J}(x;\nu +n)=\mathcal{J}(x;\nu )f_{n,\nu }(%
\frac{1}{x})-\mathcal{J}(x;\nu -1)f_{n-1,\nu +1}(\frac{1}{x}).
\label{recger}
\end{equation}%
The dual orthogonality relation of the polynomials $f_{n,\nu }(x)$ is 
\begin{equation*}
\sum_{n=0}^{\infty }\frac{B_{\nu +1}}{2\lambda _{n}(\nu +1)}f_{n,\nu +1}(%
\frac{1}{x_{\nu ,k}})f_{n,\nu +1}(\frac{1}{x_{\nu ,j}})=\frac{x_{\nu ,j}^{2}%
}{A_{j}(\nu +1)}\delta _{j,k}
\end{equation*}%
for some constants $A_{j}(\nu +1)$ and $\lambda _{n}(\nu +1)$. (for the
evaluation of these constants, as well as other parts of the argument missed
in this brief sketch, we recommend the reading of section 4 of \cite{Ism2001}%
). From (\ref{recger}) and the above analysis we can obtain the
interpolation property 
\begin{equation*}
\mathcal{J}(x_{n,\nu };\nu +k)=\frac{-\mathcal{J}(x_{n,\nu };\nu -1)}{C_{\nu
}...C_{\nu +n-1}}f_{k-1,\nu +1}(\frac{1}{x_{n,\nu }}).
\end{equation*}%
Therefore, in the language of the second section we can set 
\begin{equation*}
\mathcal{J}_{k}(t)=\sqrt{\lambda _{k}(\nu )}\frac{B_{k+\nu }}{B_{\nu }}%
\mathcal{J}(t,k)
\end{equation*}%
\begin{equation*}
r_{k}(t)=\sqrt{\lambda _{k}(\nu )}\frac{B_{k+\nu }}{B_{\nu }}f_{n,\nu }(t)
\end{equation*}%
\begin{equation*}
\lambda _{n}=\frac{-\mathcal{J}(x_{n,\nu };\nu -1)}{C_{\nu }...C_{\nu +n-1}}
\end{equation*}%
and define the kernel 
\begin{equation*}
K(x,t)=\sum_{k=0}^{\infty }u_{k}\sqrt{\lambda _{k}(\nu )}\frac{B_{k+\nu }}{%
B_{\nu }}\mathcal{J}(t,k)p_{k}(x)
\end{equation*}%
where $\left\vert u_{k}\right\vert =1$ and $\{p_{n}(x)\}$ is an arbitrary
complete orthonormal system in $L^{2}(\mu )$. As before, the kernel $K(x,t)$
can be used to define an integral transformation between two Hilbert spaces.
We could now apply the machinery of section 2 and provide a reproducing
kernel structure and a sampling theorem by means of an integral transform
with the above kernel. However, no simplification would occur on the absence
of proper addition formulas for the kernel $K(x,t)$. Choosing families of
orthogonal polynomials $f_{n,\nu }(t)$ and $p_{n}(x)$ in a way that such
addition formulas exist is the topic of the second problem in \cite{Iprob}.

Operators weighted by the Jacobi weight function \ $(1-x)^{\alpha
}(1+x)^{\beta }$\ on the interval $\left[ -1,1\right] $ can be studied by
using the Jacobi polynomials and the Wimp polynomials in a similar fashion
to what was done in section 4 of \cite[section 4]{IZ}. Reasoning as before,
we can use the resulting formulas to generalize the results in the third
section to Fourier systems with Jacobi weights, although no major
simplification seems to occur. Results very similar to those of section 3
would follow, with an extra parameter; expansions in series of $_{1}F_{1}$
replace the Neumann expansions and sampling theorems with sampling points
located at the zeros of these $_{1}F_{1}$ can also be derived.

In section 6 of \cite{IRZ}, formula (6.13) is a generalization of (\ref{qGeg}%
). This formula involves continuous $q$-analogues of the Jacobi polynomials
defined via the Askey Wilson polynomials \cite{AW} and a $q$-exponential
function with an extra variable. This more general and complicated situation
should be studied in a future work once more summation formulas are known.
Also, following Ismail%
\'{}%
s suggestion in Problem 24.2.1 of \cite{Ismbook} may lead to the discovery
of other examples of structures fitting to the setting described in this
paper.

\textbf{Aknowledgment:} I thank Hans Feichtinger for his kind hospitality
during my stay at NuHag, University of Vienna, where part of this work was
done. I%
\'{}%
m also indebted to Jose Carlos Petronilho and Fethi Bouzeffour for important
discussions on earlier versions of the manuscript.

\end{document}